\documentclass{article}
\usepackage{graphicx} 
\usepackage{amssymb}
\usepackage{amsthm}
\usepackage[cmex10]{amsmath}
\usepackage{adjustbox} 
\usepackage{breakcites} 
\usepackage{color}
\usepackage{epsfig}
\usepackage{epstopdf}
\usepackage{etoolbox}
\usepackage{enumitem}
\usepackage{float}
\usepackage{fancyhdr}
\usepackage{graphics}
\usepackage[hidelinks]{hyperref}
\usepackage{latexsym,amsfonts}
\usepackage{longtable}
\usepackage{listings}
\usepackage{lscape} 
\usepackage{lipsum}
\usepackage{multirow}             
\usepackage{pdfpages}
\usepackage[figuresright]{rotating} 
\usepackage{sectsty}
\usepackage{setspace}
\usepackage{subfigure}
\usepackage{textcomp}
\usepackage{url} 
\usepackage{wrapfig} 
\usepackage{wasysym} 
\usepackage{xcolor}
\usepackage{xcolor, soul}
\sethlcolor{yellow}
\usepackage{braket}
\usepackage{stmaryrd}
\usepackage{breqn}
\usepackage{comment}
\usepackage[margin=1in]{geometry}
\newcommand{\order}{{\mathcal O}}
\newcommand{\fluctint}[1]{\llbracket {#1} \rrbracket}
\newcommand{\avrged}[1]{\langle {#1} \rangle}
\newcommand{\avrgzero}[1]{\underline{#1}^{(0)}}

\numberwithin{equation}{section}
\newcommand\braced[1]{\left\{{#1}\right\}}
\newcommand{\imh}{i-1/2}

\title{Homogenized Equations for Isentropic Gas in a Pipe with Periodically-Varying Cross-Section}
\author{Laila S. Busaleh\thanks{{\texttt{lsbusaleh@iau.edu.sa}, Department of Mathematics, College of Science and Humanities, Jubail, Imam Abdulrahman Bin Faisal University, Saudi Arabia} \\
    Applied Mathematics and Computational Science, CEMSE Division, King Abdullah University of Science and Technology (KAUST), Thuwal, 23955-6900, Kingdom of Saudi Arabia}
\and  David I. Ketcheson\thanks{{\texttt{david.ketcheson@kaust.edu.sa}}, Applied Mathematics and Computational Science, CEMSE Division, King Abdullah University of Science and Technology (KAUST), Thuwal, 23955-6900, Kingdom of Saudi Arabia}}
\begin{document}

\maketitle
\begin{abstract}
We analyze the behavior of an isentropic gas in a narrow pipe with periodically-varying cross-sectional area.
Using multiple-scale perturbation theory, we derive homogenized effective equations, which take the form of a constant-coefficient system of evolution equations, including dispersive higher-order derivative terms.  We provide an approximate Riemann solver for the variable-cross-section isentropic gas equations, and compare numerical solutions of the original system with those of the homogenized system.  We observe that the resulting solutions take the form of solitary waves, rather than shock waves, under fairly general conditions.
\end{abstract}

\section{Introduction}
Composite materials with microstructures play a crucial and influential role in various industries. Their usage has markedly expanded in recent years, encompassing sectors such as aeronautics, aerospace, naval, and automobile manufacturing \cite{engineers}.
 Homogenization was initially developed specifically for the study of composite materials with periodic structures \cite{homo_development}. It is a valuable tool in numerous scientific and technological fields where there is a need to solve partial differential equations (PDEs) that arise in heterogeneous media. Due to material heterogeneity, these PDEs feature rapidly-oscillating coefficients. These problems are challenging to solve by analytical or computational means. From an analytical standpoint, our available tools for directly analyzing and solving PDEs with varying coefficients are quite limited. Meanwhile, numerical solution of such problems requires an extremely fine grid to accurately capture the small-scale variations. This computational requirement is costly and may render direct numerical treatment of the problem impractical. Through homogenization theory one can obtain a system with constant coefficients that approximates the original variable-coefficient system. The derivation of homogenized equations can be accomplished through various approaches, including the multiple-scale perturbation theory, Bloch wave expansion, energy estimates method, and other emerging methodologies. In this paper, we will concentrate specifically on homogenization using the multiple-scale perturbation theory.

Nonlinear wave equations can broadly be divided into two types.  Solutions of first-order nonlinear hyperbolic PDEs, such as the inviscid Burgers equation, are characterized by shock formation. In contrast, dispersive nonlinear wave PDEs,  such as the Korteweg-deVries equation (KdV) and Benjamin-Bona-Mahony equation (BBM), typically exhibit solitary wave solutions. Generally, the presence of a nonlinear term leads to steepening, resulting in shock formation. On the other hand, dispersive terms introduce oscillations that prevent shocks and give rise to solitary waves. Therefore, the long-term behavior of nonlinear waves is significantly influenced by the presence of dispersive regularizing term(s) \cite{deLunaKetch_cylindrical}.

The behavior of hyperbolic PDEs with (spatially) periodically-varying coefficients can differ significantly from that of the same PDEs with constant coefficients. For instance, the solution of the linear wave equation with periodic coefficients exhibits dispersion, even without the presence of dispersive terms in the equation. This intriguing phenomenon has been analyzed by Santosa \& and Symes \cite{santosa_symes} by employing Bloch expansions, and in other works using multiple-scale perturbation theory \cite{quezada_dispersion,ockendon2015,allaire2022crime}.
For nonlinear waves, this effective dispersion can balance with nonlinearity to generate solitary waves.
LeVeque and Yong observed this effect for nonlinear elastic waves in a periodic layered medium \cite{leveque2003}. Recent studies have explored various scenarios related to the interplay of nonlinear and effective dispersive effects  \cite{deLunaKetch_cylindrical, dreiden2010longitudinal,simpson2011coherent, 2012_ketchesonleveque_periodic, Diffractons,lombard2014numerical,chassagne2019dispersive, 2021_solitary,patil2022strongly,ketcheson2023multiscale}. 

Herein, we study the behavior of an isentropic gas propagating within a pipe with a periodically-varying cross-section. Similar to the studies just mentioned, we observe that effective dispersion from the periodic structure balances with nonlinear effects and leads to the formation of solitary waves in typical solutions.


While the periodic variations of the pipe's cross-section introduce complexity, we can capture the dominant wave behavior through a simplified model. This model uses effective parameters to capture the average effects of the microstructure, simplifying the analysis without losing the essential dynamics of wave propagation.

The main goal of this work is to construct an effective medium that represents the essential characteristics of the original microstructured medium.
To accomplish this, we develop a constant-coefficient (homogenized) approximation, which effectively describes the behavior of long-wavelength solutions for a system that describes the scenario of a gas propagating in a narrow pipe with a periodically-varying cross-sectional area. This approximation method draws upon the approach developed by Yong and co-authors \cite{yong2002, leveque2003}, and more recently was applied to the shallow water waves over periodically-varying bathymetry \cite{ketcheson2023multiscale,ketcheson2025dispersive},
as well as the compressible Euler equations \cite{ketcheson2025solitary}.
The constant-coefficient equations we derive have dispersive terms, demonstrating that waves in this context exhibit effective dispersion and allowing us to quantify the dispersive effects.
While the formal process employed herein is
very similar to that in the works just mentioned, the
analysis here is more technical due to the combined presence of variable coefficients and a variable non-conservative term (see Equation \eqref{isothermal-pipe}).  The resulting homogenized equations are, correspondingly, more complicated than in any of the previous applications.

The outline of this paper is as follows. In Section 2, we describe the model problem and perform the multiple-scale analysis which leads to the effective equations.  The linear dispersion relation of the homogenized equations is determined in section 3. In section 4, we solve an initial value problem and compare the results of the original problem with those of the effective model.

Code to reproduce all simulation results in this paper can be found online\footnote{\url{https://github.com/laila-16/Homogenization_of_Isentropic_Gas_Equations}}.
\section{Model Equations}
In this paper we consider the one-dimensional flow of an isentropic gas in a pipe, modeled by the first-order hyperbolic system \cite{lefloch2003riemann}:
\begin{subequations} \label{isothermal-pipe}
\begin{align}
   (a\rho)_t + (a \rho u)_x & = 0, \\
    (a \rho u)_t + (a\rho u^2 + a P(\rho))_x &= P(\rho)a_x.
\end{align}
\end{subequations}
Here $\rho(x,t)$ and $u(x,t)$ represent the density and velocity, respectively, and $a(x)> 0$ is the cross-sectional area that varies periodically with period $\delta$: $$a(x+\delta) = a(x).$$ Finally, $P(\rho)$ represents the pressure, and it is given by$$
    P(\rho) = \kappa \rho^\gamma,  ~~ \text{where} ~~ 1<\gamma<5/3.
$$\begin{figure}[htb]
    \center
    \includegraphics[width=5.5in]{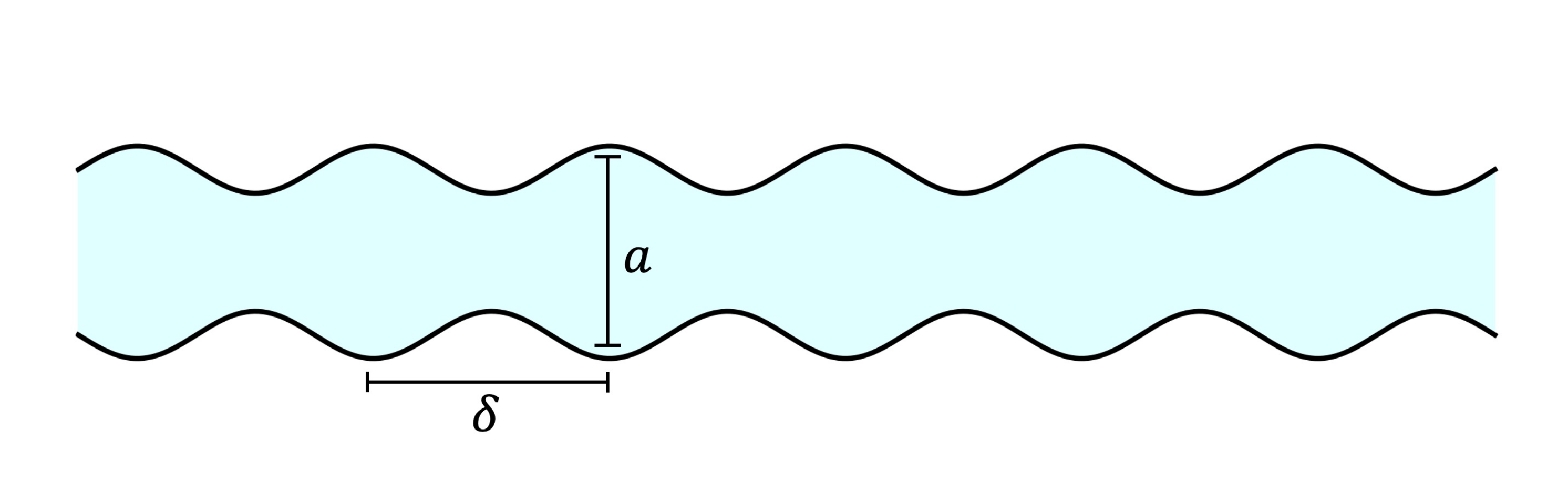}
    \caption{Flow in a pipe of periodically-varying cross-section.\label{pipe-flow}}
\end{figure}
Herein we take $\gamma= 1.4$, though this specific value is not essential to our analysis or results.
In this system, the velocity $u$ varies rapidly for any non-trivial near-steady flow as $a(x)$ varies. This rapid variation is not conducive to the homogenization process. Therefore, we introduce the flux $q=a\rho u$ which provides a better starting point for the homogenization, since the flux  and the density can vary slowly even for solutions in which $a(x)$ varies rapidly. Thus, we rewrite \eqref{isothermal-pipe} in terms of $q$ as
\begin{subequations} \label{iso-pipe-reformulated}
\begin{align}
    a\rho_t + q_x & = 0, \\
    q_t + (q^2/(a\rho))_x + a P'(\rho)\rho_x & = 0.
\end{align}
\end{subequations}
In the analysis, it is convenient to assume that $a(x)$ is a smooth function, although we will see in the numerical examples that the homogenized equations we derive below are effective even when $a(x)$ is discontinuous.

\subsection{Ansatz}
Now to find the homogenized equations, we introduce a microscopic length variable that we refer to as the fast spatial scale:
$$y=x/\delta,$$
where again $\delta$ is the spatial period of the variation in $a(x)$.
Formally, $x$ and $y$ are assumed to be independent, so we  transform the spatial derivative as:
\begin{align} \label{eq:dtrans}
\frac{\partial}{\partial x} \to \frac{\partial}{\partial x} + \delta^{-1} \frac{\partial}{\partial y},
\end{align}
and we treat $a$ as a function of $y$ only.
Applying \eqref{eq:dtrans} to \eqref{iso-pipe-reformulated} yields:
\begin{subequations}  
\label{transformed-system}
\begin{align}
\label{first_eq_trans-sys}
 a\rho_t + q_x+ \delta^{-1} q_y&=0\\
 \label{sec_eq_trans-sys}
q_t +\left(\frac{q^2}{a\rho}\right)_x + \delta^{-1}\left(\frac{q^2}{a\rho}\right)_y +a P'(\rho)\rho_x+\delta^{-1}aP'(\rho)\rho_y&=0 
\end{align}
\end{subequations}
We assume that $q$ and $\rho$ can be expressed in terms of power series in $\delta$:
\begin{align}
\label{asymptotic_expansion}
q(x,y,t;\delta)= {q}^{(0)}(x,t)+\delta q^{(1)}(x,y,t)+\delta^2 q^{(2)}(x,y,t)+\delta^3 q^{(3)}(x,y,t)+... \\ \nonumber
\rho(x,y,t;\delta)= {\rho}^{(0)}(x,t)+\delta \rho^{(1)}(x,y,t)+\delta^2 \rho^{(2)}(x,y,t)+\delta^3 \rho^{(3)}(x,y,t)+...
\end{align}
The functions $q^{(i)}$ and $\rho^{(i)}$ are assumed to be sufficiently regular and periodic with respect to $y$, with period $\delta$, in order to avoid the appearance of secular terms in what follows. 
For simplicity, and without loss of generality,
we take $\delta=1$.
We will also make use of the Taylor series expansion:
\begin{align} \label{dpdrho_expansion}
    P'(\rho)=P'(\rho^{(0)})+P''(\rho^{(0)})\left(\delta \rho^{(1)}+\delta^2 \rho^{(2)}+...\right)+\frac{1}{2}P'''(\rho^{(0)})\left(\delta \rho^{(1)}+\delta^2 \rho^{(2)}+...\right)^2 +...
\end{align}
We substitute \eqref{asymptotic_expansion} and \eqref{dpdrho_expansion} into \eqref{transformed-system} and collect terms proportional to each power of $\delta$.
To save space, we omit the lengthy result here and simply present the collected terms proportional to each power of $\delta$, below.

\subsection{Averaging Operators}
Before continuing the asymptotic analysis we introduce the following linear operators that
we will use in the averaging processes to identify terms that are independent of $y$.
\begin{itemize}
        \item $\avrged{b(y)} = \int^{1}_{0} b(y) dy.$  This averaging operator gives the average value of a function over one period. The result of this operator is a function independent of $y$.  
        \item $\{b\}(y) = b(y) - \avrged{b(y)}$. This operator generates the part of the function that has zero average (the fluctuating part). It returns a function that depends on $y$.
        \item $\fluctint{b}(y) = \left\{\int_0^y \braced{b(\xi)}\, d\xi\right\}.$ This operator gives the integral of the fluctuating part of a function, where the constant of integration is chosen such that the average of the integral of the fluctuating part is zero. This operator returns a function that depends on $y$. 
\end{itemize}
Some properties of the averaging operators are introduced in \cite{ketcheson2023multiscale,yong2002}.

\subsection{Asymptotic Analysis}

With the necessary definitions introduced, we can now proceed with the process of term collection. We will demonstrate the derivation of the homogenized equations for the first few orders of $\delta$. 
As the expressions are lengthy, we will omit the explicit function dependencies, but we will
underline quantities that are independent of the fast variable $y$.
\subsubsection{$\order(\delta^{-1})$}
From \eqref{first_eq_trans-sys} there are no terms proportional to $\delta^{-1}$. From \eqref{sec_eq_trans-sys} there is only one term, yielding 
\begin{subequations}
\label{avg-1_order}
    \begin{align}
     -\frac{(\underline{q}^{(0)})^2 a_y}{a^2\underline{\rho}^{(0)}}=0. 
    \end{align}
\end{subequations}
Since in general $a_y\neq 0$, we conclude that $\underline{q}^{(0)} =0$.

\subsubsection{$\order(\delta^0)$}
After setting $\underline{q}^{(0)} =0$, the $\order(\delta^0)$ equations are
\begin{subequations}
\label{zero-order-system}
    \begin{align} 
     a \underline{\rho}^{(0)}_t +  q^{(1)}_y =0 , \\
     \label{zero-order-system-secondeq}
     a P'(\underline{\rho}^{(0)})\rho^{(1)}_y + a P'(\underline{\rho}^{(0)})\underline{\rho}^{(0)} _x = 0.
    \end{align}
\end{subequations}
We assume $a P'(\underline{\rho}^{(0)}) \neq 0$, so \eqref{zero-order-system-secondeq} reduces to 
$$ \rho^{(1)}_y + \underline{\rho}^{(0)} _x = 0.$$
By solving for  $q^{(1)}_y$ and $\rho^{(1)}_y$ we get:
\begin{subequations}
\label{zero-order-system-y}
    \begin{align}
    \label{first-zero-ordereq}
      q^{(1)}_y &= - a\underline{\rho}^{(0)}_t , \\
    \label{second-zero-ordereq}
    \rho^{(1)}_y &= -\underline{\rho}^{(0)} _x .
    \end{align}
\end{subequations}
Since $q^{(i)}, \rho^{(i)}$ are periodic in $y$ for $i>0$, integrating \eqref{first-zero-ordereq} over one period with respect to $y$ gives
$$0=\int_0^1 - a \underline{\rho}^{(0)}_t dy =  - \underline{\rho}^{(0)}_t\int_0^1 a dy = -\avrged{a} \underline{\rho}^{(0)}_t,$$ 
which implies $\underline{\rho}^{(0)}_t=0$. Similarly, we integrate \eqref{second-zero-ordereq} over one period with respect to $y$, obtaining $\underline{\rho}^{(0)} _x=0$.
Therefore, $\underline{\rho}^{(0)}$ is constant: $\underline{\rho}^{(0)}(x,t)=\underline{\rho}^{(0)}$. 
We conclude furthermore that $q^{(1)}(x,y,t)= \underline{q}^{(1)}(x,t)$ and $\rho^{(1)}(x,y,t)= \underline{\rho}^{(1)}(x,t)$.
Hence the asymptotic expansion \eqref{asymptotic_expansion} reduces to: 
\begin{align*}
    q(x,y,t;\delta)&= \delta \underline{q}^{(1)}(x,t)+\delta^2 q^{(2)}(x,y,t)+\delta^3 q^{(3)}(x,y,t)+...\\
    \rho(x,y,t;\delta)&= \underline{\rho}^{(0)}+\delta\underline{\rho}^{(1)}(x,t)+\delta^2 \rho^{(2)}(x,y,t)+\delta^3 \rho^{(3)}(x,y,t)+ ...,
\end{align*}
and the corresponding solutions are $\order(\delta)$
perturbations from the steady state with density $\rho^{(0)}$ and velocity zero.

\subsubsection{$\order(\delta^{1})$}
The $\order(\delta^1)$ equations are
\begin{subequations}
    \begin{align}
    0&= a\underline{\rho}^{(1)}_t + \underline{q}^{(1)}_x + q^{(2)} _y ,\\[11pt]
    0&=-\frac{a^{-2}a_y}{\underline{\rho}^{(0)}}(\underline{q}^{(1)})^2+\underline{q}^{(1)}_t + aP'(\underline{\rho}^{(0)})\left(\underline{\rho}^{(1)}_x +\rho^{(2)}_y\right).
    \end{align}
\end{subequations}
Solving for $\rho^{(2)}_y$ and $q^{(2)}_y$ we get :
\begin{subequations}
\label{order-one-equations}
    \begin{align}
     -q^{(2)} _y &= a\underline{\rho}^{(1)}_t + \underline{q}^{(1)}_x ,\\[11pt]
     -\rho^{(2)}_y &= \frac{-a^{-3}a_y}{\underline{\rho}^{(0)}P'(\underline{\rho}^{(0)})}(\underline{q}^{(1)})^2 +\frac{a^{-1}}{P'(\underline{\rho}^{(0)})}\underline{q}^{(1)}_t + \underline{\rho}^{(1)}_x,
    \end{align}
\end{subequations}To remove the $y$ dependence and avoid secular terms we apply the averaging operator $\avrged{.}$. We obtain: 
\begin{subequations}
\label{order-one-equations-avg}
    \begin{align}
     \order(\delta^{1})&= \avrged{a}\underline{\rho}^{(1)}_t + \underline{q}^{(1)}_x,  \\[11pt]
    \order(\delta^{1})&= \frac{\avrged{a^{-1}}}{P'(\underline{\rho}^{(0)})} \underline{q}^{(1)}_t+\underline{\rho}^{(1)}_x.
    \end{align}
\end{subequations}
The average of the term $-\frac{a^{-3}a_y(\underline{q}^{(1)})^2}{\underline{\rho}^{(0)}P'(\underline{\rho}^{(0)})}$ vanishes because $a(y)$ is periodic in $y$ and hence $\int_{\delta} a^{-3} a_y dy =0  $. Notice that \eqref{order-one-equations-avg} is just the linearized and averaged form of the original system \eqref{iso-pipe-reformulated}. 
We define the averaged variables 
\begin{align}
\label{reduced_expansion}
    \underline{q}&= \underline{q}^{(1)} + \delta\underline{q}^{(2)}+ \delta^2\underline{q}^{(3)}+ ...\\ \nonumber
    \underline{\rho}&= \underline{\rho}^{(1)} + \delta \underline{\rho}^{(2)} + \delta^2 \underline{\rho}^{(3)} + ...
\end{align}
Hence the $\order(\delta^{1})$ averaged equations are 
\begin{subequations} 
\label{first-order-averaged}
\begin{align}
   \underline{\rho}_t & = -\frac{\underline{q}_x}{\avrged{a}} + \order(\delta)  ,\\
   \underline{q}_t &= -\frac{P'(\underline{\rho}^{(0)})}{\avrged{a^{-1}}}\underline{\rho}_x + \order(\delta) . 
\end{align}
\end{subequations}
In order to find $q^{(2)}$ and $\rho^{(2)}$ we subtract
\eqref{order-one-equations} from \eqref{order-one-equations-avg} and integrate with respect to $y$. This yields:
\begin{subequations}
\label{rho_q_2vals}
    \begin{align}
     q^{(2)}&=-\fluctint{a} \underline{\rho}^{(1)}_t + \underline{q}^{(2)},\\[11pt]
     \rho^{(2)}&= \frac{\fluctint{a^{-3}a_y}}{P'(\underline{\rho}^{(0)}) \underline{\rho}^{(0)}} (\underline{q}^{(1)})^2- \frac{\fluctint{a^{-1}}}{P'(\underline{\rho}^{(0)})} \underline{q}^{(1)}_t +\underline{\rho}^{(2)},
    \end{align}
\end{subequations}
where $\underline{q}^{(2)}(x,t)$ and $\underline{\rho}^{(2)}(x,t)$ are constants of integration in terms of $y$.
\subsubsection{$\order(\delta^2)$}
The $\order(\delta^2)$ equations are
\begin{subequations}
\begin{align}
    0&= a\rho^{(2)}_t + q^{(2)}_x + q^{(3)}_y ,
    \\[11pt]
    0&=  q^{(2)}_t + a^{-2}a_y \left(-\frac{2q^{(1)}q^{(2)}}{\rho^{(0)}} + \frac{(q^{(1)})^2 \rho^{(1)}}{(\rho^{(0)})^2}\right) + 2 a^{-1} \left(\frac{q^{(1)} q^{(2)}_y + q^{(1)}q^{(1)}_x}{\rho^{(0)}}\right)\nonumber \\[11pt] 
    &\ \  
    +a P''(\rho^{(0)})\rho^{(1)}\left(\rho^{(2)}_y + \rho^{(1)}_x\right)+ aP'(\rho^{(0)})\left(\rho^{(3)}_y + \rho^{(2)}_x\right).
\end{align}    
\end{subequations}
Next we replace ${\rho}^{(1)}, {\rho}^{(2)}, {q}^{(1)}, {q}^{(2)}$ and their derivatives by their values from  \eqref{rho_q_2vals} and \eqref{zero-order-system-y}.
By solving the resulting equations for $q^{(3)}_y$ and $\rho^{(3)}_y$ we obtain: 
\begin{subequations}
\label{order-two-equations}
\begin{align}
    -q^{(3)}_y&= a\left(\frac{2\fluctint{a^{-3}a_y}}{P'(\underline{\rho}^{(0)}) \underline{\rho}^{(0)}}\underline{q}^{(1)}\underline{q}^{(1)}_t - \frac{\fluctint{a^{-1}}}{P'(\underline{\rho}^{(0)})}\underline{q}^{(1)}_{tt} +\underline{\rho}^{(2)}_t \right) -\fluctint{a} \underline{\rho}^{(1)}_{xt} + \underline{q}^{(2)}_x,
    \\[11pt]
    -\rho^{(3)}_y&=  \frac{-a^{-1}\fluctint{a} }{P'(\underline{\rho}^{(0)})}\underline{\rho}^{(1)}_{tt} + \frac{a^{-1}}{P'(\underline{\rho}^{(0)})} \underline{q}^{(2)}_t+ a^{-3}a_y \left(\frac{2\fluctint{a} \underline{q}^{(1)}\underline{\rho}^{(1)}_t -2\underline{q}^{(1)}\underline{q}^{(2)}}{P'(\underline{\rho}^{(0)})\rho^{(0)}} + \frac{(q^{(1)})^2 \rho^{(1)}}{P'(\underline{\rho}^{(0)})(\rho^{(0)})^2}\right) \nonumber \\[11pt] 
    &\ \  - \frac{2a^{-1}}{P'(\underline{\rho}^{(0)})\rho^{(0)}} \underline{q}^{(1)} \underline{\rho}^{(1)}_t +P''(\underline{\rho}^{(0)})\rho^{(1)}\left(\frac{a^{-3}a_y}{\underline{\rho}^{(0)}(P'(\underline{\rho}^{(0)})^2}(\underline{q}^{(1)})^2-\frac{a^{-1}}{(P'(\underline{\rho}^{(0)}))^2} \underline{q}^{(1)}_t \right)\nonumber \\[11pt] 
    &\ \ + \frac{2\fluctint{a^{-3}a_y}}{P'(\underline{\rho}^{(0)}) \underline{\rho}^{(0)}} \underline{q}^{(1)}\underline{q}^{(1)}_x - \frac{\fluctint{a^{-1}}}{P'(\underline{\rho}^{(0)})} \underline{q}^{(1)}_{xt}+\underline{\rho}^{(2)}_x.
\end{align}    
\end{subequations}
Applying the averaging operator $\avrged{.}$ to the system above in order to remove the $y$ dependence yields the following averaged system:
\begin{subequations}  
\label{order-two-equations-avg}
\begin{align}   
\label{order-two-equations-firstavg}
 \nonumber \\[8pt] 
  \order(\delta^2)&=  \frac{2\avrged{a\fluctint{a^{-3}a_y}}}{P'(\underline{\rho}^{(0)}) \underline{\rho}^{(0)}}\underline{q}^{(1)}\underline{q}^{(1)}_t - \frac{\avrged{a\fluctint{a^{-1}}}}{P'(\underline{\rho}^{(0)})} \underline{q}^{(1)}_{tt}+\avrged{a}\underline{\rho}^{(2)}_t + \underline{q}^{(2)}_x,
    \\[13pt]
    \label{order-two-equations-secondavg}  \order(\delta^2)&=  -\frac{ \avrged{a^{-1}\fluctint{a}}}{P'(\underline{\rho}^{(0)})}\underline{\rho}^{(1)}_{tt} + \frac{\avrged{a^{-1}}}{P'(\underline{\rho}^{(0)})}\underline{q}^{(2)}_t + \frac{2 \avrged{a^{-3}a_y \fluctint{a}}}{P'(\underline{\rho}^{(0)})\underline{\rho}^{(0)}} \underline{q}^{(1)}\underline{\rho}^{(1)}_t  \nonumber \\[11pt] 
    &\ \  - \frac{2\avrged{a^{-1}} }{P'(\underline{\rho}^{(0)})\underline{\rho}^{(0)}} \underline{q}^{(1)} \underline{\rho}^{(1)}_t -\frac{\avrged{a^{-1}} P''(\underline{\rho}^{(0)})}{(P'(\underline{\rho}^{(0)}))^2} \underline{\rho}^{(1)}\underline{q}^{(1)}_t +\underline{\rho}^{(2)}_x .
\end{align}        
\end{subequations}
Next, in order to write these equations as an evolution system, we 
use \eqref{first-order-averaged} and \eqref{order-two-equations-avg} to replace the second-order time derivatives and mixed derivatives by spatial derivatives only. The $\order(\delta^2)$ averaged equations are
\begin{subequations} 
\label{secondorderaveraged}
\begin{align}
\label{secondorderaveraged1}
   \underline{\rho}_t & = -\frac{\underline{q}_x}{\avrged{a}}  + \delta \left(- \frac{2\avrged{a\fluctint{a^{-3}a_y}} }{\avrged{a^{-1}} \avrgzero{\rho} \avrged{a}} \underline{q}\,\underline{\rho}_x -\frac{\avrged{a\fluctint{a^{-1}}} }{\avrged{a^{-1}} \avrged{a}^2}\underline{q}_{xx}\right) + \order(\delta^{2}),\\
    \nonumber \\[11pt]
    \label{secondorderaveraged2}
   \underline{q}_t &= -\frac{P'}{\avrged{a^{-1}}}\underline{\rho}_x +\delta \Biggr[\left(\frac{-2}{\avrgzero{\rho} \avrged{a}}+\frac{2\avrged{a^{-3}a_y\fluctint{a}}}{\avrged{a^{-1}}\avrgzero{\rho}\avrged{a}} \right)\underline{q}\,\underline{q}_x -\frac{P''}{\avrged{a^{-1}}}\underline{\rho}\,\underline{\rho}_x +\frac{\avrged{a^{-1}\fluctint{a}} P'}{\avrged{a^{-1}}^2 \avrged{a}}\underline{\rho}_{xx}\Biggl]+ \order(\delta^{2}) .
\end{align}
\end{subequations}
The second-derivative terms in the above equations turn out to be
dispersive in nature.  However, their coefficients turn out to be quite small, or to vanish identically for many choices of $a(x)$, and so in typical scenarios (like those considered in Section \ref{sec:numerical}) the dominant dispersion will come from the next order terms.  
In order to be able to find the dominant dispersive terms at the next order, we subtract
\eqref{order-two-equations} from \eqref{order-two-equations-avg} and integrate with respect to $y$ to get :

\begin{subequations} \label{ttt-form}
\begin{align}
    q^{(3)}&= -\frac{2\fluctint{a\fluctint{a^{-3}a_y}}}{P'(\underline{\rho}^{(0)}) \underline{\rho}^{(0)}}\underline{q}^{(1)}\underline{q}^{(1)}_t +\frac{\fluctint{a\fluctint{a^{-1}}}}{P'(\underline{\rho}^{(0)})} \underline{q}^{(1)}_{tt}-\fluctint{a}\underline{\rho}^{(2)}_t  +\fluctint{\fluctint{a}} \underline{\rho}^{(1)}_{xt} + \underline{q}^{(3)},
    \\[11pt]
    \rho^{(3)}&=  \frac{\fluctint{a^{-1}\fluctint{a}} }{P'(\underline{\rho}^{(0)})} \underline{\rho}^{(1)}_{tt} - \frac{\fluctint{a^{-1}}}{P'(\underline{\rho}^{(0)})}\underline{q}^{(2)}_t + \left( \frac{2\fluctint{a^{-1} }}{P'(\underline{\rho}^{(0)})\underline{\rho}^{(0)}}- \frac{2\fluctint{a^{-3}a_y\fluctint{a}} }{P'(\underline{\rho}^{(0)})\underline{\rho}^{(0)}} \right)\underline{q}^{(1)}\underline{\rho}^{(1)}_t  \nonumber \\[11pt] 
    &\ \    -\left(\frac{\fluctint{a^{-3}a_y}}{P'(\underline{\rho}^{(0)})(\avrgzero{\rho})^2} + \frac{\fluctint{a^{-3}a_y}P''(\underline{\rho}^{(0)})}{\underline{\rho}^{(0)}(P'(\underline{\rho}^{(0)})^2} \right) (\underline{q}^{(1)})^2\underline{\rho}^{(1)} + \frac{\fluctint{a^{-1}}P''(\underline{\rho}^{(0)})}{(P'(\avrgzero{\rho}))^2} \underline{q}^{(1)}_t\underline{\rho}^{(1)}\nonumber \\[11pt] 
    &\ \ + \frac{2\fluctint{a^{-3}a_y}}{P'(\underline{\rho}^{(0)})\underline{\rho}^{(0)}}\underline{q}^{(1)}\underline{q}^{(2)} - \frac{2\fluctint{\fluctint{a^{-3}a_y}}}{P'(\underline{\rho}^{(0)}) \underline{\rho}^{(0)}}\underline{q}^{(1)}\underline{q}^{(1)}_x + \frac{\fluctint{\fluctint{a^{-1}}}}{P'(\underline{\rho}^{(0)})} \underline{q}^{(1)}_{xt} +\underline{\rho}^{(3)}.
\end{align}    
\end{subequations}

\subsubsection{Governing Equations} \label{sec:governing}
To derive the $\order(\delta^{3})$ averaged equations, we follow a procedure similar to the one used for lower-order terms. Specifically, we obtain the $\order(\delta^{3})$ governing equations by combining $\delta$ times \eqref{first-order-averaged}, $\delta^{2}$ times \eqref{secondorderaveraged}, and $\delta^{3}$ times the $\order(\delta^{3})$ equations (the detailed expressions are omitted here due to their complexity). We then simplify the expression and combine terms using \eqref{reduced_expansion}.

The resulting equations contain terms with high-order time derivatives.  A linear dispersion analysis of the equations reveals that they are unstable.
This is a known issue with certain homogenization techniques like that employed herein; see \cite{allaire2022crime} for a general discussion and \cite{leveque2003,ketcheson2023multiscale} for similar examples.
 What is recommended in \cite{allaire2022crime},
 and also employed in \cite{leveque2003}, is to exchange all high-order time derivatives (i.e., all except those appearing in the linear evolution term for each equation) for space derivatives, by using equality of mixed partial derivatives and keeping the same formal order of accuracy in $\delta$.
 This approach leads to a system that is also linearly unstable, but in a weaker sense -- it exhibits instability only for high-wavenumber modes, and consequently it is possible to use it for numerical simulations (as is done in \cite{leveque2003}) as long as the mesh is not too fine.

This yields
\begin{subequations} 
\label{thirdorderaveraged}
\begin{align}
   \label{thirdorderaveraged1}
   \underline{\rho}_t & = -\frac{\underline{q}_x}{\avrged{a}}  +\delta \left(-\frac{C_{1} }{\avrged{a^{-1}} \avrged{a}^{2}} \underline{q}_{xx}+ \frac{2C_{13}}{\avrged{a^{-1}} \avrged{a}} \underline{q}\,\underline{\rho}_{x} \right)  \\[12pt]& \ \ \nonumber 
   + \delta^2 \Biggl[\left(\frac{4C_{13}^2 + 4C_{13}C_{14}}{\avrged{a^{-1}}P'\avrged{a}^2}
   - \frac{C_3}{P'\avrged{a}^2}- \frac{C_{13}P''}{(P')^2 \avrged{a}^{2}}\right)\underline{q}^{2} \underline{q}_{x} + \left(\frac{C_{9}}{\avrged{a^{-1}} \avrged{a}^{3}}-\frac{C_{2}}{\avrged{a^{-1}}\avrged{a}^{2}}\right)\underline{q}_{xxx} \\[12pt]& \ \ \nonumber  
   -\frac{2C_{3}}{\avrged{a^{-1}}\avrged{a}}\underline{q}\,\underline{\rho}\,\underline{\rho}_{x} + \left(\frac{2C_{8}}{\avrged{a^{-1}}\avrged{a}^{2}} - \frac{2C_{4}}{\avrged{a^{-1}}\avrged{a}}\right)\underline{q}_{x} \underline{\rho}_{x}  +\left( -\frac{2C_{13} C_{1}}{\avrged{a^{-1}}^{2} \avrged{a}^{2}} +\frac{2 C_{8}}{\avrged{a^{-1}} \avrged{a}^{2}} -\frac{2C_{4}}{\avrged{a^{-1}}\avrged{a}}\right)\underline{q}\,\underline{\rho}_{xx}\Biggr]+\order(\delta^{3}),
   \\[23pt]
   \label{thirdorderaveraged2}
   \underline{q}_t &= -\frac{P'}{\avrged{a^{-1}}}\underline{\rho}_{x} +\delta \Biggl[\frac{-2C_{13}-2C_{14}}{\avrged{a^{-1}}\avrged{a}} \underline{q}\,\underline{q}_{x} -\frac{P''}{\avrged{a^{-1}}}\underline{\rho}\,\underline{\rho}_{x} +\frac{C_1 P'}{\avrged{a^{-1}}^{2} \avrged{a}}\underline{\rho}_{xx}\Biggr]  \\[12pt]& \ \ \nonumber +\delta^{2} \Biggl[\frac{2C_{10}+2C_{3}}{\avrged{a^{-1}} \avrged{a}}\underline{q}\,\underline{\rho}\,\underline{q}_{x}  +\left(\frac{ 2C_{1} C_{13}}{\avrged{a^{-1}}^{2} \avrged{a}^{2}}-\frac{C_{8}}{\avrged{a^{-1}}\avrged{a}^{2}}\right)\underline{q}_{x}^{2} \\[12pt]& \ \ \nonumber +\left(-\frac{2C_{1}C_{13}}{\avrged{a^{-1}}^{2} \avrged{a}^2}+\frac{4C_{7}}{\avrged{a^{-1}}^{2} \avrged{a} }-\frac{2C_{4}}{\avrged{a^{-1}}\avrged{a}} \right)\underline{q}\,\underline{q}_{xx} + \left(\frac{-3C_{5}+4C_{6}+C_{12}}{\avrged{a^{-1}}^{2}}+\frac{4C_{13}^2 + 4C_{13}C_{14}}{\avrged{a^{-1}}^{2} \avrged{a}}\right)\underline{q}^{2}\underline{\rho}_x \\[12pt]& \ \ \nonumber +\left(\frac{2C_{7}P'}{ \avrged{a^{-1}}^{3}}+\frac{C_{1}P''}{\avrged{a^{-1}}^{2} \avrged{a}}-\frac{2 C_{1} C_{13} P'}{\avrged{a^{-1}}^{3} \avrged{a}}\right)\underline{\rho}_x^{2} 
   +\frac{C_{1}P''}{\avrged{a^{-1}}^{2} \avrged{a}}\underline{\rho}\,\underline{\rho}_{xx} +\left(\frac{C_{11}P'}{\avrged{a^{-1}}^{3} \avrged{a}}-\frac{C_{2} P'}{\avrged{a^{-1}}^{2} \avrged{a}}\right)\underline{\rho}_{xxx}\Biggr] + \order(\delta^{3}),
\end{align}
\end{subequations}
where the coefficients are provided in Appendix \ref{appA}.
We remark that for a wide class of periodic functions $a(x)$ (specifically, those that are translation-even in the terminology of \cite[Appendix~A]{ketcheson2023multiscale}) we have $C_1=C_7=0$, so that several of the coefficients above vanish.  In the examples in Section \ref{sec:numerical} we also choose such functions $a(x)$, and in some of the following analysis we will consider this special, slightly simpler case.

Similar to the other applications mentioned above, these equations are unstable for sufficiently high wavenumbers (see our analysis below).
In order to obtain a fully stable system of equations,
we follow what is done
in \cite{ketcheson2023multiscale} and was suggested by observations dating back to Whitham \cite[Section~13.11]{whitham2011linear}.  Namely, we
exchange some time derivatives for space derivatives by using
equality of mixed partial derivatives, but we
keep exactly one $t$-derivative in each higher-order linear term, thus obtaining a system that is linearly stable for all wavenumbers.  The resulting system can be
written (introducing symbols for the lengthy coefficients, which are also defined Appendix \ref{appA}):
\begin{subequations} 
\label{BBM}
\begin{align}
   \underline{\rho}_t & = -\frac{1}{\avrged{a}}\underline{q}_x  +\delta \alpha_3 \underline{q}\,\underline{\rho}_x + \delta^2 \left(\alpha_4\underline{q}^2 \underline{q}_x -\avrged{a} \alpha_{5} \underline{\rho}_{xxt} +\alpha_6 \underline{q}\,\underline{\rho}\,\underline{\rho}_x + \alpha_7 \underline{q}_x \underline{\rho}_x +\alpha_8 \underline{q}\,\underline{\rho}_{xx}\right),\\[12pt]
    \underline{q}_t &= -\frac{P'}{\avrged{a^{-1}}} \underline{\rho}_x + \delta\left(\beta_2 \underline{q}\,\underline{q}_x +\beta_3 \underline{\rho}\,\underline{\rho}_x \right) +\delta^2 \Bigl( \beta_5 \underline{q}\,\underline{\rho}\,\underline{q}_x 
   +\beta_6 \underline{q}_x^2 +\beta_7 \underline{q}\,\underline{q}_{xx} +\beta_{8}\underline{q}^2\underline{\rho}_x\\[12pt]& \ \ \nonumber +\beta_{9}\underline{\rho}_x^2-\frac{\avrged{a^{-1}}}{P'}\beta_{11}\underline{q}_{xxt}\Bigr).
\end{align}
\end{subequations}

We conclude this session with a remark regarding the ansatz \eqref{asymptotic_expansion}.  It is possible to instead suppose \emph{a priori} that the waves of interest are quadratically small with respect to $\delta$.  This leads to a similar set of equations except that the power of $\delta$ multiplying each term is different.  In particular, with such an ansatz the quadratic nonlinear terms and the third-derivative terms appear at the same order as one another.

\section{Linear Dispersion Relation} \label{section: dispersion_relation}
To quantify the effective dispersion due to the
periodic variation in $a(x)$, we study the dispersion
relation for small-amplitude solutions of the homogenized
equations \eqref{thirdorderaveraged}.
To find the dispersion relation we start by linearizing the system around the constant state $(\rho^{(0)}, 0)$. We let 
\begin{align}
\label{dispersion_sub}
    \underline{\rho} &= \rho^{(0)} + \epsilon \overline{\rho}(x,t),
    \\ \nonumber
    \underline{q} &= \epsilon \overline{q}(x,t),
\end{align}
where $\epsilon \ll1 $. Moreover, for simplicity, we introduce shorthand for the coefficients in \eqref{thirdorderaveraged}, which are detailed in Appendix \ref{appA}. Thus, \eqref{thirdorderaveraged} can be written compactly as
\begin{subequations} 
\label{govrn}
\begin{align}
  \label{govrn1}
   \underline{\rho}_t & = -\frac{1}{\avrged{a}}\underline{q}_x  + \delta \alpha_3 \underline{q}\,\underline{\rho}_x + \delta^2 \left( \alpha_4\underline{q}^2 \underline{q}_x + \alpha_5 \underline{q}_{xxx} +\alpha_6 \underline{q}\,\underline{\rho}\,\underline{\rho}_x + \alpha_7 \underline{q}_x \,\underline{\rho}_x +\alpha_8 \underline{q}\,\underline{\rho}_{xx}\right),\\[12pt]
   \label{govrn2}
    \underline{q}_t &= -\frac{P'}{\avrged{a^{-1}}} \underline{\rho}_x +\delta \left(\beta_2 \underline{q}\,\underline{q}_x +\beta_3 \underline{\rho}\,\underline{\rho}_x\right) + \delta^2 \Bigl(\beta_5 \underline{q}\,\underline{\rho}\,\underline{q}_x 
   +\beta_6 \underline{q}_x^2 +\beta_7 \underline{q}\,\underline{q}_{xx} +\beta_{8}\underline{q}^2\underline{\rho}_x +\beta_{9}\underline{\rho}_x^2\\[12pt]& \ \ \nonumber+\beta_{11}\underline{\rho}_{xxx}\Bigr).
\end{align}
\end{subequations}
Now we substitute \eqref{dispersion_sub} into \eqref{govrn} and collect terms of first order in $\epsilon$, neglecting terms of order $\order{(\epsilon^2)}$ (and higher). This results in:
\begin{subequations} 
\label{linearized_system}
\begin{align}
   \overline{\rho}_t & = -\frac{1}{\avrged{a}}\overline{q}_x + \delta^2 \alpha_5 \overline{q}_{xxx},\\[12pt]
   \overline{q}_t &= -\frac{P'}{\avrged{a^{-1}}} \overline{\rho}_x +\delta \beta_3\rho^{(0)}\overline{\rho}_x +\delta^2 \beta_{11}\overline{\rho}_{xxx}.
\end{align}
\end{subequations}
The obtained system is linear with constant coefficients, so we can represent its solution by a linear combination of Fourier modes. Let:
\begin{align}
\label{FT}
\overline{\rho}=\overline{q}= e^{i(k x - \omega t)},
\end{align}
where $k$ denotes the wave number and $\omega$ is the frequency. This form is referred to as a mode of the linear evolution equation. We plug \eqref{FT} into \eqref{linearized_system} to get: 
\begin{subequations}
\label{fourier_sub}
\begin{align}
   -i\omega \overline{\rho}&=  -\frac{1}{\avrged{a}}  ik\overline{q} - ik^3 \delta^2 \alpha_5  \overline{q},
    \\[15pt] 
    -i\omega\overline{q}&=  -\frac{P'}{\avrged{a^{-1}}}ik \overline{\rho}+ ik \delta\beta_3\rho^{(0)} \overline{\rho} - ik^3 \delta^2 \beta_{11} \overline{\rho},
\end{align}
\end{subequations}
which has the form
\begin{align*}
\begin{bmatrix}
\omega  & \hfill  -\frac{1}{\avrged{a}} k 
   - k^3 \delta^2 \alpha_5  \\
  -\frac{P'}{\avrged{a^{-1}}}k +k \delta \beta_3\rho^{(0)} - k^3 \delta^2 \beta_{11} \hfill &  \omega
\end{bmatrix}
\begin{bmatrix}
\overline{\rho}\\
\overline{q}
\end{bmatrix} := 
    M\begin{bmatrix}
\overline{\rho}\\
\overline{q}
\end{bmatrix}=\begin{bmatrix}
0\\
0
\end{bmatrix}.
\end{align*}
The dispersion relation is obtained by setting $\det M = 0$, which yields
\begin{align} \label{xxx_dispersion-relation}
    \omega(k) = \pm |k|\sqrt{\left(-\frac{1}{\avrged{a}} -k^2 \delta^2 \alpha_5\right)\left(-\frac{P'}{\avrged{a^{-1}}} + \delta \beta_3 \rho^{(0)} - k^2 \delta^2 \beta_{11} \right)}.
\end{align}
For the problem to be stable we require the dispersion relation $\omega (k)$ to be real for all (real) wave numbers $k$ \cite{Deconinck}. 
However, $\omega(k)$ in \eqref{xxx_dispersion-relation} is in fact complex for large enough wavenumbers.
This is an expected issue for equations obtained through the process used here (see \cite{allaire2022crime,ketcheson2023multiscale}, where the same phenomenon is observed).
While it's possible to deal with this issue in numerical simulations through filtering, a more robust approach is as follows.

Recall that in Section \ref{sec:governing} we modified the original homogenized equations \eqref{ttt-form} by exchanging $t$-derivatives for $x$-derivatives, for instance replacing a $q_{ttt}$ term by one involving $\rho_{xxx}$.  It turns out to be advantageous to modify this process and
retain one $t$-derivative in each of the highest-order dispersive terms.  The better properties of
these mixed-derivative dispersive terms is known in the context of the BBM and KdV equations.

Equivalently, working backward from
the equations \eqref{thirdorderaveraged}, we make the following replacements using relationships derived from the first-order averaged equations \eqref{first-order-averaged}, as follows:
\begin{subequations} 
\begin{align}
   \underline{\rho}_{txx} & = -\frac{\underline{q}_{xxx}}{\avrged{a}}  \iff \underline{q}_{xxx} = -\avrged{a} \underline{\rho}_{txx} ,\\
   \underline{q}_{txx} &= -\frac{P'}{\avrged{a^{-1}}}\underline{\rho}_{xxx} \iff  \underline{\rho}_{xxx} = -\frac{\avrged{a^{-1}}}{P'} \underline{q}_{txx}
\end{align}
\end{subequations}
Then the homogenized equations become:
\begin{subequations} 
\label{BBM_truncated}
\begin{align}
   \overline{\rho}_t & = -\frac{1}{\avrged{a}}\overline{q}_x -\avrged{a}\delta^2\alpha_{5} \overline{\rho}_{xxt},\\[12pt]
   \overline{q}_t &= -\frac{P'}{\avrged{a^{-1}}} \overline{\rho}_x +\delta \beta_3 \rho^{(0)} \overline{\rho}_x -\frac{\avrged{a^{-1}}}{P'}\delta^2 \beta_{11} \overline{q}_{xxt}.
\end{align}
\end{subequations}
Therefore, the dispersion relation is 
\begin{align} \label{disprel-xxt}
    \omega(k) = \pm |k|\sqrt{\frac{-\frac{1}{\avrged{a}} \left(-\frac{P'}{\avrged{a^{-1}}} + \delta \beta_3 \rho^{(0)} \right)}{\left(1-\avrged{a}\delta^2\alpha_{5}k^2 \right)\left(1-\frac{\avrged{a^{-1}}}{P'}\delta^2\beta_{11}k^2 \right)}}
\end{align}
The dispersion relation is stable when \(\mathrm{Im}(\omega(k)) = 0\). Since the product of \(-\frac{1}{\avrged{a}}\) and \(\left(-\frac{P'}{\avrged{a^{-1}}} + \delta \beta_3 \rho^{(0)}\right)\) is always positive for any \(\rho^{(0)} \geq 0\), we can conclude that the condition for stability is satisfied. Specifically, \(-\frac{1}{\avrged{a}}\), \(-\frac{P'}{\avrged{a^{-1}}}\), and \(\beta_3\) depend only on \(\avrged{a}\) and the derivatives of the pressure, all of which are positive. For the denominator of the expression, $-\avrged{a}\alpha_{5}$ and $-\frac{\avrged{a^{-1}}}{P'}\beta_{11}$ are always positive for any periodic function, further ensuring stability. 
For the numerator, we can show that both \(C_9\) and \(C_{11}\) are negative. Starting with \(C_9\), we have:
\[
C_9 = \avrged{a \fluctint{a^{-1} \fluctint{a}}} = - \avrged{a^{-1} \fluctint{a}^2},
\]
the averaged quantity is clearly positive, so multiplying it by a negative sign will always yield a negative result. Therefore, $C_9 $ is negative, and similarly, it can be shown that $C_{11}$ will remain negative for all values. Additionally, we observe that \(C_2 > 0\) for the cross-sections we tested, and this condition is fulfilled for the cross-section chosen in this paper. Thus, the stability conditions are met for the selected cross-section and parameters used in this paper.

\begin{figure}[htb]
    \centering
    \includegraphics[width=\textwidth, trim=0cm 0.5cm 0cm 0.3 cm, clip]{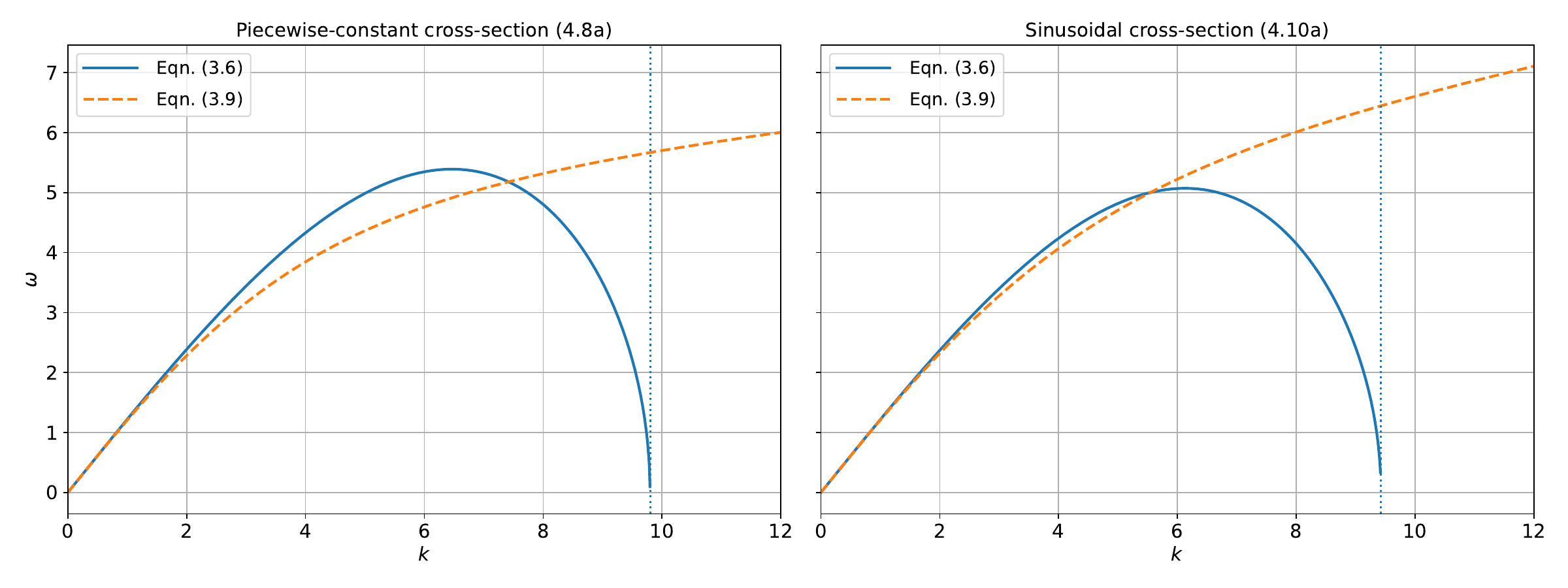}
    \caption{Dispersion relations with piecewise-constant (left) and sinusoidal (right) cross-sections. Solid lines correspond to KdV-type dispersion~\eqref{xxx_dispersion-relation}, and dashed lines to BBM-type dispersion~\eqref{disprel-xxt}. 
    Note that for \eqref{xxx_dispersion-relation}, we truncate the plot at the point where $\omega(k)$ becomes complex.
    }
    \label{dispersion-comparison}
\end{figure}

In Figure~\ref{dispersion-comparison}, we compare the dispersion relations for the specific cross-sections under consideration. We observe that at large wave numbers, the dispersion relation given by~\eqref{xxx_dispersion-relation} grows significantly, which can impact the stability of the numerical solution. To analyze the agreement between the two models, we perform a Taylor expansion of both dispersion relations.

For both dispersion relations (3.6) and (3.9), the Taylor expansion gives  
$$\omega(k) = \sqrt{-\frac{1}{\avrged{a}}  (-\frac{P'}{\avrged{a^{-1}}} + \delta\beta_3 \rho_0)k^2} + \mathcal{O}(k^3),$$  
since the $\mathcal{O}(k^2)$ terms vanish in each case.  
Thus, the two dispersion relations coincide up to the included order, and their difference satisfies  
$\omega_{\text{(3.6)}}(k) - \omega_{\text{(3.9)}}(k) = \mathcal{O}(k^3)$.

\section{Numerical Comparison} \label{sec:numerical}
In this section we propose numerical discretizations for
both the variable-coefficient isentropic gas system \eqref{isothermal-pipe}
and the homogenized equations \eqref{BBM}.  We then compare numerical solutions of the two systems in order to understand the accuracy of the homogenized approximation.

\subsection{Numerical Discretization of Homogenized System}

We discretize the full nonlinear homogenized equations \eqref{BBM}, rewriting them first the form
\begin{subequations}
\label{pseudo}    
\begin{align} 
    \label{rhot}\overline{\rho}_t & =  (1- \alpha_{5b} \partial_x^2)^{-1}\mathcal{F}_1(\overline{\rho},\overline{q}), \\
    \label{qt}\overline{q}_t & = (1- \beta_{11b} \partial_x^2)^{-1} \mathcal{F}_2(\overline{\rho},\overline{q}),
\end{align}
\end{subequations}
Here, $\mathcal{F}_1(\overline{\rho}, \overline{q})$ and $\mathcal{F}_2(\overline{\rho}, \overline{q})$ are the right-hand sides of the first and second equations in \eqref{BBM}, respectively, excluding the $xxt$-derivative terms.
We solve this system using a Fourier pseudospectral collocation method in space and a
third-order strong stability preserving Runge-Kutta method \cite{shu-osher} in time.

\subsection{Numerical Discretization of Isentropic Gas Equations}
To derive an approximate Riemann solver for the isentropic gas equation, we rewrite \eqref{isothermal-pipe} in terms of the momentum \(\mathtt{m} = u \rho\) and density \(\rho\) as
\begin{subequations} \label{iso-reformulated}
\begin{align}
    a\rho_t + (a\mathtt{m})_x & = 0, \\
    a\mathtt{m}_t + \left( \frac{a\mathtt{m}^2}{\rho} + aP(\rho) \right)_x  & = P(\rho) a_x.
\end{align}
\end{subequations}
The system \eqref{iso-reformulated} has the general form of a conservation law with a source term:
\begin{align}
\label{capacity}
    a(x) V_t + f(V,x)_x = \Psi(x),
\end{align}
where
\begin{align}
V & =
\begin{bmatrix}
   \rho \\ \mathtt{m} 
\end{bmatrix}, &
f(V,x) & = \begin{bmatrix}
    a(x)\mathtt{m} \\ \frac{a(x)\mathtt{m}^2}{\rho} + a(x)P(\rho)
\end{bmatrix}, &
\Psi= 
 \begin{bmatrix}
   0 \\
    a_x P(\rho)
\end{bmatrix}.
\end{align}
The flux $f(V,x)$ depends explicitly on both \(x\) and the vector of conserved quantities \(V\). 

In \eqref{iso-reformulated}, the conserved quantity is \(a(x)V\), and \(a(x)\) is sometimes referred to as a \emph{capacity function} \cite{FVM}. 
The system also has a spatially-varying flux function and source term. 
We use an f-wave Riemann solver and capacity-form
differencing \cite{fwave,FVM} in order to solve this system
in a way that handles these complications and allows us to accurately compute solutions that are close to a steady state (i.e., the discretization is well-balanced).

The finite volume method approximates
the cell averages
$$V^{n}_i = \frac{1}{\Delta x} \int_{x_{i - 1/2}}^{x_{i + 1/2}} v(x, t_n) \,dx $$ where $\Delta x = x_{i + 1/2} - x_{i - 1/2} $ is the length of the cell.  At each cell interface we solve a Riemann problem with the piece-wise initial data
$$v(x,0) = \begin{cases} V_{i-1} & x < x_{i-1/2} ,  \\
   V_i & x > x_{i-1/2}.
    \end{cases}$$
The solution is updated as:
$$V^{n+1}_{i}= V^n_i - \frac{\Delta t}{ \Delta x}\left(\mathcal{A}^+ \Delta V_{i-1/2} + \mathcal{A}^- \Delta V_{i+1/2}\right) $$\\ 
where the fluctuations are given by: 
\begin{subequations} \label{fluctuations}
\begin{align}
    \mathcal{A}^{+} \Delta V_{i-1/2} & =\sum_{p:s^p_\imh > 0} \mathcal{Z}_{i-1/2}^p,  \\
    \mathcal{A}^{-} \Delta V_{i+1/2} & =\sum_{p:s^p_\imh < 0} \mathcal{Z}_{i+1/2}^p,
\end{align}
\end{subequations}
 where $p= 1, 2.$ The terms $\mathcal{A}^{-} \Delta V_{i+1/2}$ and $\mathcal{A}^{+} \Delta V_{i-1/2}$ represent the contributions to the cell average $V_i$ due to right-going waves from $x_{i-1/2}$ and left-going waves from 
 $x_{i+1/2}$, respectively.

It remains only to prescribe the numerical approximation of the waves.
 As the system is genuinely nonlinear, the exact solution of the Riemann problem involves two waves, each of which may be a shock or a rarefaction.
We approximate this solution by two propagating discontinuities (herein referred to as waves) that are determined as follows.
We approximate the flux Jacobian at \(x_\imh\) by
\begin{equation}
\label{jacobi}
A_\imh = \begin{bmatrix}
    0 & \tilde{a}  \\
    -\frac{\tilde{a} \tilde{\mathtt{m}} ^2}{\tilde{\rho} ^2} + \tilde{a} P'(\tilde{\rho}) & \frac{2\tilde{a} \tilde{\mathtt{m}} }{\tilde{\rho} }
\end{bmatrix}.
\end{equation}
Here, we use the averaged values 
$$\tilde{a} \equiv \frac{a_i+a_{i-1}}{2},\quad \tilde{\rho} \equiv \frac{\rho_i+\rho_{i-1}}{2}\text{\quad and \quad}\tilde{\mathtt{m}} \equiv \frac{\mathtt{m}_i+\mathtt{m}_{i-1}}{2}.$$
We discretize the source term at the cell interfaces as
\[
\Psi_{i-1/2} = \frac{1}{2}\left(P_{i}(\rho) + P_{i-1}(\rho)\right) \frac{a_i - a_{i-1}}{\Delta x}.
\]
To determine the waves, we decompose 
the flux difference after subtracting the source
term: 
\[
f_i(V_i) - f_{i-1}(V_{i-1}) - \Delta x \Psi_{i-1/2} = \sum_{p=1}^{2} \beta_{i-1/2}^p r_{i-1/2}^p .
\]
Here $r_{i-1/2}^p$ are the right eigenvectors of
$A_\imh$:
\[
R_\imh = 
\begin{bmatrix}
    \frac{\tilde{\rho} }{\tilde{\mathtt{m}} - \tilde{\rho}\sqrt{P'(\tilde{\rho})}} & \frac{\tilde{\rho}}{\tilde{\mathtt{m}} + \tilde{\rho}\sqrt{P'(\tilde{\rho})}} \\
    1 & 1
\end{bmatrix}
 = \begin{bmatrix} r^1_\imh & r^2_\imh \end{bmatrix}.
\]
In order to determine which wave(s) contribute to each
fluctuation in \eqref{fluctuations}, we also need approximate speeds,
for which we use the corresponding eigenvalues of $A_\imh$:
\[
\lambda^1_\imh =\text{min}\left( a_l \frac{{\mathtt{m}}_l }{\rho_l } - a_l \sqrt{P'(\rho_l)}, a_r \frac{{\mathtt{m}}_r }{\rho_r } - a_r \sqrt{P'(\rho_r)}\right), \] \[\lambda^2_\imh= \text{max}\left( a_l \frac{{\mathtt{m}}_l }{\rho_l } + a_l \sqrt{P'(\rho_l)}, a_r \frac{{\mathtt{m}}_r }{\rho_r } + a_r \sqrt{P'(\rho_r)}\right),
\]

The waves and speeds required to compute
the fluctuations \eqref{fluctuations} are thus
\begin{align}
    \mathcal{Z}^p_{i-1/2} & = \beta^p_{i-1/2}r_{i-1/2}^p & s^p= \lambda^p_\imh.
\end{align}

\subsection{Numerical tests}
\subsubsection{Piecewise-Constant Cross-section}
\begin{figure}[h]
    \centering
    \includegraphics[width=\textwidth, trim=0cm 4.9cm 0cm 0.5cm, clip]{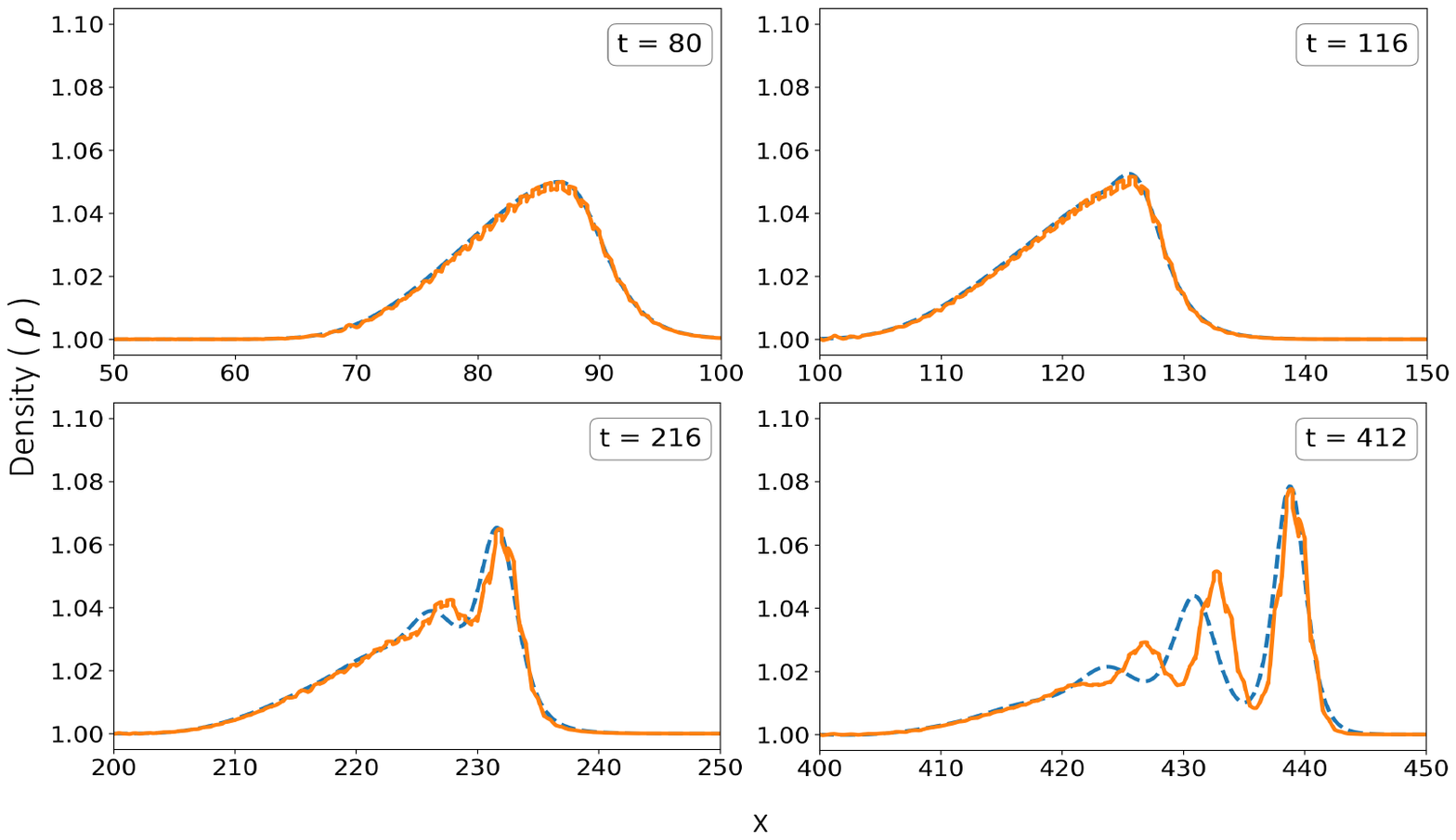}
    \caption{First experiment ($\delta=0.1$): comparison of the homogenized system (dashed line) and the variable-coefficient system (solid line) solutions, with the cross-section and initial data given by \eqref{scenario_a}.}
    \label{fig:scenario_a_delta}
\end{figure}

\begin{figure}[h]
    \centering
    \includegraphics[width=\textwidth, trim=0cm 4.9cm 0cm 0.5cm, clip]{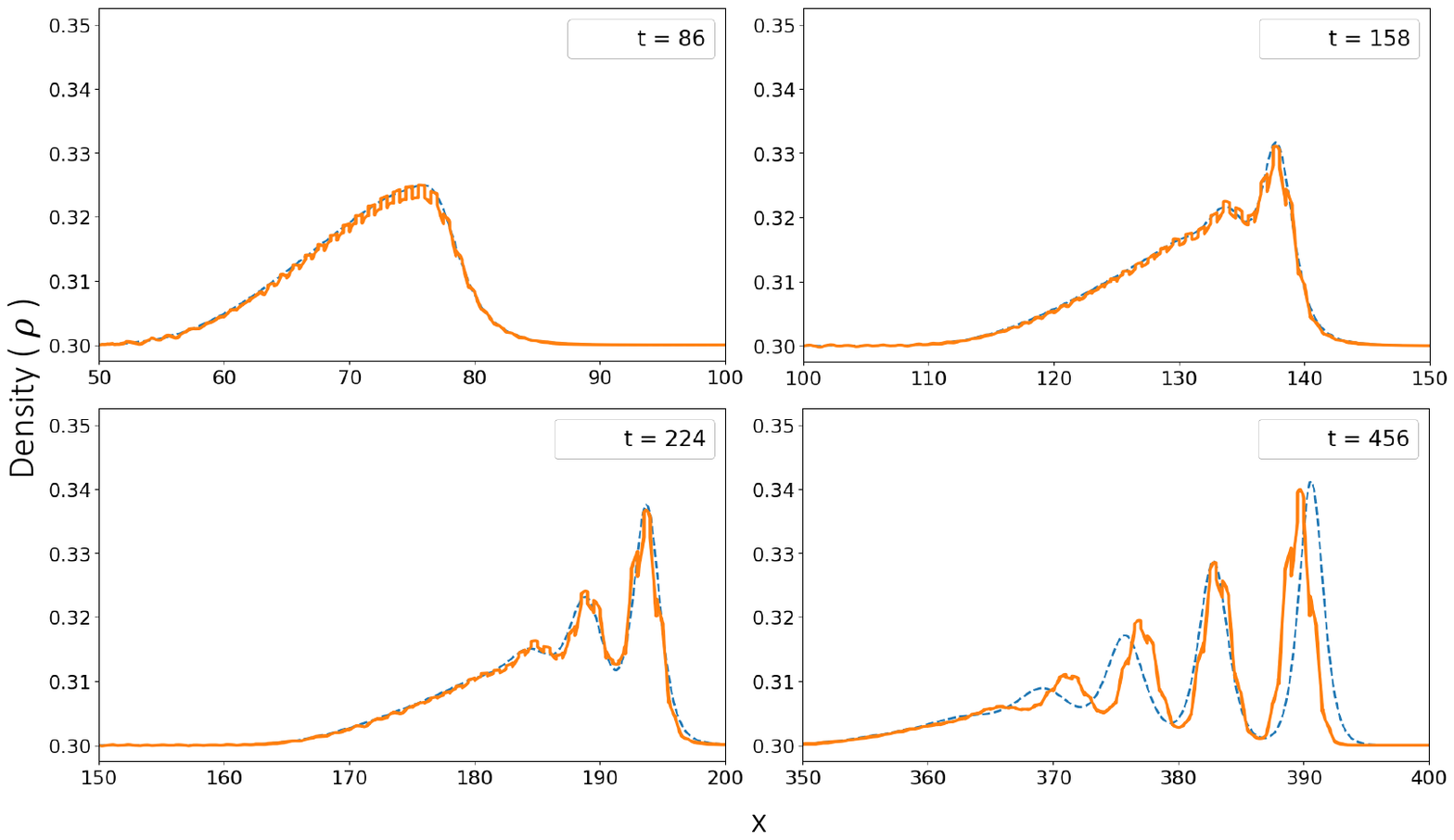}
    \caption{Second experiment: comparison of the homogenized system (dashed line) and the variable-coefficient system (solid line) solutions, with the cross-section and initial data given by \eqref{scenario_a}.}
    \label{fig:scenario_a}
\end{figure}

In this section, we consider flow in a pipe with piecewise-constant cross-section, where the cross-section and initial data are defined as follows. We conduct two different experiments corresponding to two choices of initial conditions.

\begin{subequations} \label{scenario_a}
\begin{align}
    a(x) & = \begin{cases} 
    \frac{1}{4}, & 0 \leq x - \lfloor x \rfloor < \frac{1}{2}, \\
    \frac{3}{4}, & \frac{1}{2} \leq x - \lfloor x \rfloor < 1, 
    \end{cases} \\
    \rho(x,0) & = 1 + \delta\, e^{-(x/8)^2}, \\
    q(x,0) & = 0,
\end{align}
\end{subequations}
Here, $0<\delta\ll 1$ is the Gaussian amplitude; we take $\delta=0.1$. In the first experiment, the Gaussian perturbation is superimposed on a constant background of $\rho=1$, so that the initial pulse divides into left- and right-going parts, each of which evolve into a train of traveling waves.
In the second experiment, we instead take
\begin{equation}
    \rho(x,0) = 0.3 + \tfrac{1}{20} e^{-(x/8)^2}, \qquad q(x,0) = 0,
\end{equation}
which produces a similar splitting of the initial pulse, but with a smaller background state. 
Figures \ref{fig:scenario_a_delta} and \ref{fig:scenario_a} show the right-going part for the first and second experiments, respectively. In both cases, we compare the pseudospectral solution of the homogenized approximation with the finite volume solution of the variable cross-section isentropic gas system \eqref{iso-pipe-reformulated}. We observe that the homogenized solution provides a good approximation at early times but, as expected, becomes less accurate as time progresses. 

\subsubsection{Smoothly-varying Cross-section}

\begin{figure}[h]
    \centering
    \includegraphics[width=\textwidth, trim=0cm 4.7cm 0cm 0.5cm, clip]{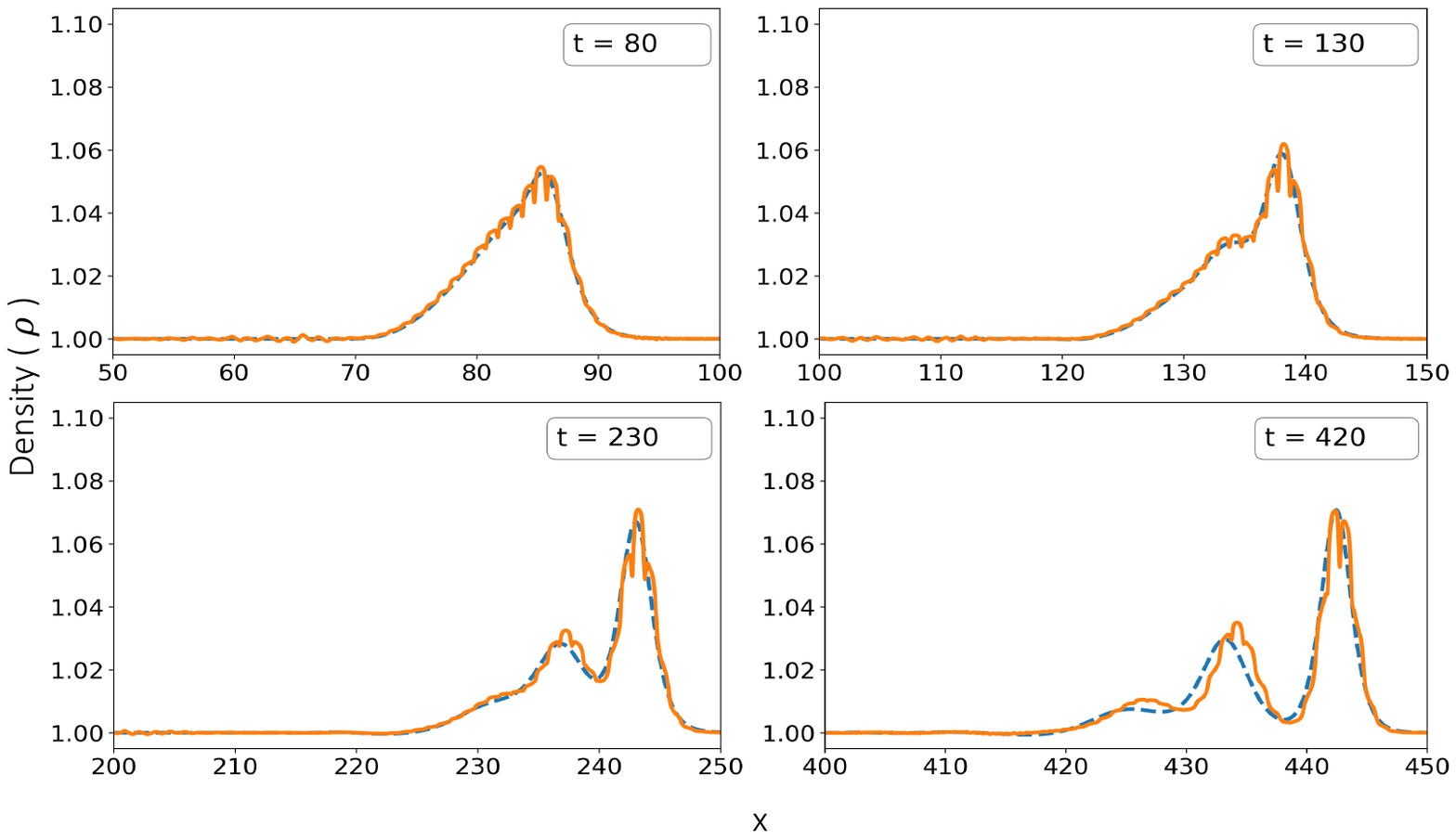}
    \caption{First experiment ($\delta=0.1$): comparison of the homogenized system (dashed line) and the variable coefficients system (solid line) solutions, with the cross-section and initial data given by \eqref{scenario_b}.}
    \label{fig:scenario_b_delta}
\end{figure}

\begin{figure}[h]
    \centering
    \includegraphics[width=\textwidth, trim=0cm 4.7cm 0cm 0.5cm, clip]{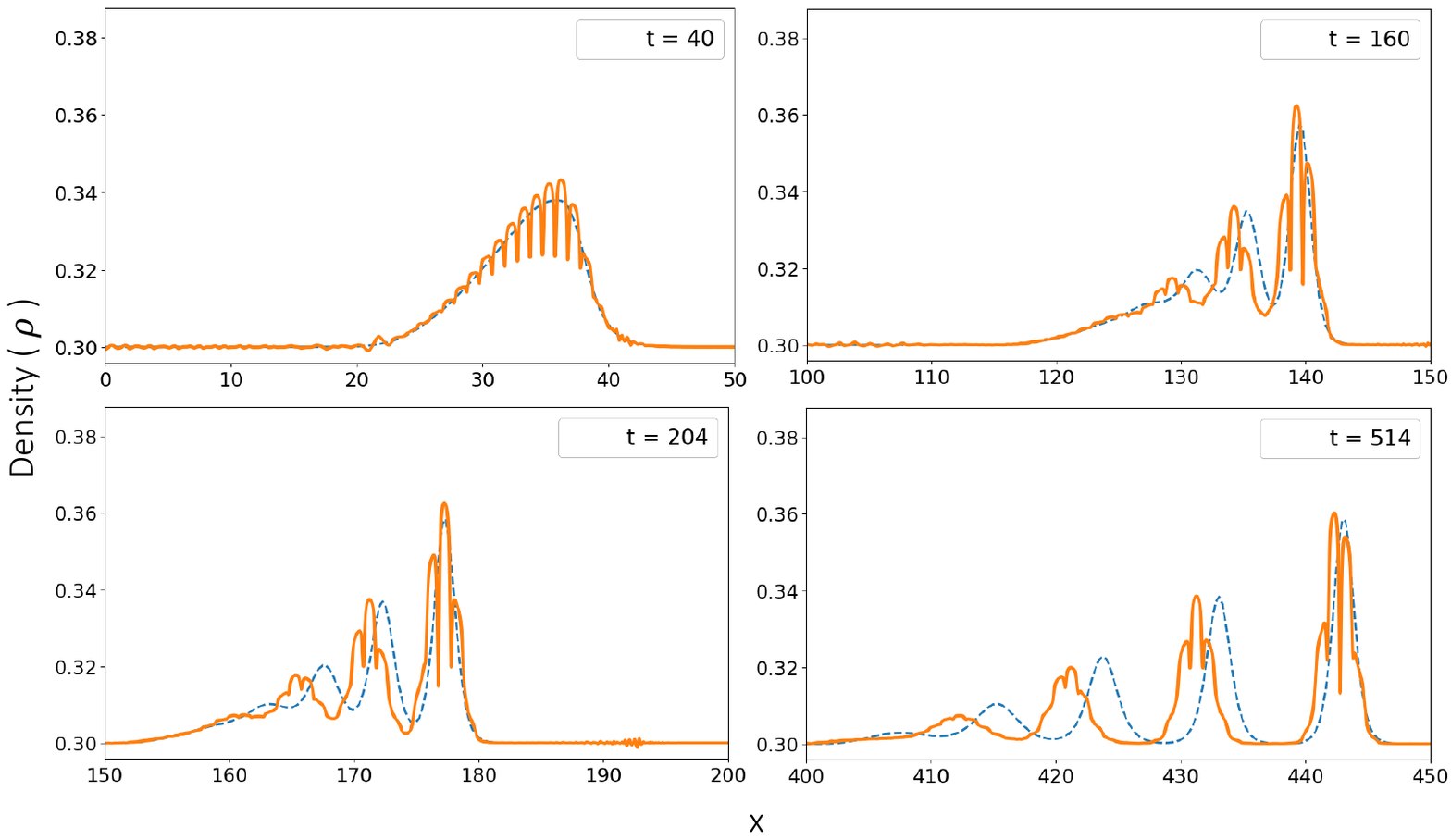}
    \caption{Second experiment: comparison of the homogenized system (dashed line) and the variable coefficients system (solid line) solutions, with the cross-section and initial data given by \eqref{scenario_b}.}
    \label{fig:scenario_b}
\end{figure}

Next, we consider a scenario with a smoothly-varying cross-section. We conduct two different experiments corresponding to two choices of initial conditions.

\begin{subequations} \label{scenario_b}
\begin{align}
    a(x) & = \frac{3}{5} + \frac{2}{5} \sin(2 \pi x), \\
    \rho(x,0) & = 1 + \delta\, e^{-(x/5)^2}, \\
    q(x,0) & = 0,
\end{align}
\end{subequations}
Here, we use the same value of $\delta$ as in the piecewise-constant cross-section. In the first experiment, the Gaussian perturbation is superimposed on a constant background of $\rho=1$.  In the second experiment, we instead take
\begin{equation}
    \rho(x,0) = 0.3 + \tfrac{1}{12} e^{-(x/5)^2}, \qquad q(x,0) = 0,
\end{equation}
which adds the Gaussian perturbation to a smaller background state.  Figures \ref{fig:scenario_b_delta} and \ref{fig:scenario_b} show the results for the first and second experiments, respectively. In both cases, we compare the pseudospectral solution of the homogenized approximation with the finite volume solution of the variable cross-section isentropic gas system \eqref{iso-pipe-reformulated}. Similar to the piecewise constant case discussed earlier, the accuracy of the homogenized approximation decreases over time.

In both scenarios, the chosen cross-sections and initial data are designed to promote the formation of solitary waves. For significantly larger initial data or for cross-sections with smaller variations relative to the average fluid depth, wave breaking occurs, causing the homogenized model to inadequately capture the full solution dynamics \cite{ketcheson2023multiscale}.

\section{Conclusions}
In this study, we derived effective equations with constant coefficients (homogenized equations) to approximate the first-order hyperbolic system of isentropic gas equations. Our analysis revealed the formation of solitary waves during the propagation of nonlinear waves in a heterogeneous medium. The presence of dispersive effects, induced by the medium's heterogeneity, leads to the decomposition of shocks into a series of solitary traveling waves.

We demonstrated this phenomenon theoretically by deriving the linear dispersion relation for the homogenized equations, which allowed us to conclude that the phase velocity is wave-number dependent. 

The derivation of the homogenized system was rigorously verified using \emph{Mathematica}. We employed the pseudospectral method for spatial discretization to solve the homogenized equations. For the original system, we used the capacity-form differencing approach and an approximate Riemann solver based on the f-wave strategy, particularly suited for solutions close to equilibrium between the flux and source terms.

By comparing the results obtained from the homogenized system with those from the variable coefficients system, we observed significant agreement, thereby validating the accuracy and effectiveness of this approximation.

\section*{Acknowledgment}
This work was supported by funding from King Abdullah University of Science and Technology.

\appendix
\section{Coefficients of the Homogenized Equations}\label{appA}

In this section, we present the coefficients for the homogenized equations. The coefficients of \eqref{thirdorderaveraged} are given by:

\begin{flalign*}
    \vspace{5em}\hspace{1.75em}C_1 ~&= \avrged{a^{-1}\fluctint{a}}, &\\
    \vspace{5em}\hspace{1.75em}C_2 ~&= \avrged{a^{-1}\fluctint{\fluctint{a}}}, &\\
    \vspace{5em}\hspace{1.75em}C_3 ~&= \frac{1}{ (\avrgzero{\rho})^2}\avrged{a\fluctint{a_ya^{-3}}}  = \frac{-1}{2 (\avrgzero{\rho})^2}\left(\avrged{a^{-1}} - \avrged{a}\avrged{a^{-2}} \right), &\\
    \vspace{5em}\hspace{1.75em}C_4 ~&=  \frac{1}{ (\avrgzero{\rho})}\avrged{a^{-3}a_y\fluctint{\fluctint{a}}}
    = \frac{-1}{2\avrgzero{\rho}}\avrged{a\fluctint{a^{-2}}}, &\\
    \vspace{5em}\hspace{1.75em}C_5 ~&=  \frac{1}{ (\avrgzero{\rho})^2}\avrged{a^{-1}\fluctint{a^{-3}a_y}}
   = \frac{-1}{2(\avrgzero{\rho})^2}\left(\avrged{a^{-3}} - \avrged{a^{-1}}\avrged{a^{-2}} \right), &\\
    \vspace{5em}\hspace{1.75em}C_6 ~&= \frac{1}{ (\avrgzero{\rho})^2}\avrged{a^{-3}a_y\fluctint{a\fluctint{a^{-3}a_y}}}
   = \frac{-1}{4(\avrgzero{\rho})^2}\left(\avrged{a^{-3}} + \avrged{a}\avrged{a^{-2}}^2 - 2\avrged{a^{-1}}\avrged{a^{-2}} \right), &\\
    \vspace{5em}\hspace{1.75em}C_7 ~&= \frac{1}{ (\avrgzero{\rho})}\avrged{a^{-3}a_y\fluctint{a\fluctint{a^{-1}}}}
    = \frac{1}{2\avrgzero{\rho}}\avrged{a^{-2}} C_1, &\\
    \vspace{5em}\hspace{1.75em}C_8 ~&= \frac{1}{ (\avrgzero{\rho})}\avrged{a\fluctint{a^{-3}a_y\fluctint{a}}}, &\\
    \vspace{5em}\hspace{1.75em}C_9 ~&= \avrged{a\fluctint{a^{-1}\fluctint{a}}}, &\\
    \vspace{5em}\hspace{1.75em}C_{10} &=\frac{1}{ (\avrgzero{\rho})^2}\avrged{a^{-1}}, &\\
    \vspace{5em}\hspace{1.75em}C_{11} &= \avrged{a^{-1}\fluctint{a\fluctint{a^{-1}}}}, &\\
    \vspace{5em}\hspace{1.75em}C_{12} &= \frac{1}{ (\avrgzero{\rho})^2}\avrged{a^{-3}}, &\\
    \vspace{5em}\hspace{1.75em}C_{13} &= \frac{1}{ (\avrgzero{\rho})}C_3, &\\
    \vspace{5em}\hspace{1.75em}C_{14} &=\frac{1}{ (\avrgzero{\rho})}\avrged{a^{-1}}. &
     &
\end{flalign*}
We make some simplifications using the identities
\begin{align*}
\avrged{f\fluctint{g}} & = - \avrged{g\fluctint{f}} \\
\avrged{f\fluctint{\fluctint{g}}} & = \avrged{g\fluctint{\fluctint{f}}}\\
\fluctint{f^{-n}f'} & = \frac{\avrged{f^{1-n}} - f^{1-n}}{n-1}.
\end{align*}

Finally, the coefficients of \eqref{govrn1} are:
\begin{align*}
    \alpha_1 &= \frac{-1}{\avrged{a}}, ~~~~~~~~~~~~~~~~~\alpha_2 = \frac{-C_{1} }{\avrged{a^{-1}} \avrged{a}^{2}} , ~~~~~~~~~~~~~~~~~~~~~~~~~ \alpha_3  = \frac{2C_{13}}{\avrged{a^{-1}} \avrged{a}} ,\\[12pt]
    \alpha_4 & =  \frac{4C_{13}^2 + 4C_{13}C_{14}}{\avrged{a^{-1}}P'\avrged{a}^2}
   - \frac{C_3}{P'\avrged{a}^2}- \frac{C_{13}P''}{(P')^2 \avrged{a}^{2}}, ~~~~~~~~~\alpha_5 =  \frac{C_{9}}{\avrged{a^{-1}} \avrged{a}^3}-\frac{C_{2}}{\avrged{a^{-1}}\avrged{a}^2},  \\[12pt]\alpha_6 &= \frac{-2C_{3}}{\avrged{a^{-1}}\avrged{a}},~~~~~~~~~ \alpha_7 = \frac{2C_{8}}{\avrged{a^{-1}}\avrged{a}^{2}} - \frac{2C_{4}}{\avrged{a^{-1}}\avrged{a}}, ~~~~~
    \alpha_8 =   -\frac{2C_{13} C_{1}}{\avrged{a^{-1}}^{2} \avrged{a}^{2}} +\frac{2 C_{8}}{\avrged{a^{-1}} \avrged{a}^{2}} -\frac{2C_{4}}{\avrged{a^{-1}}\avrged{a}},
\end{align*}
and the coefficients of \eqref{govrn2} are:
\begin{align*}
     \beta_1& =\frac{-P'}{\avrged{a^{-1}}},~~~~~~~  \beta _2= \frac{-2C_{13}-2C_{14}}{\avrged{a^{-1}}\avrged{a}}, ~~~~~~ \beta _3=  \frac{-P''}{\avrged{a^{-1}}}, \\[12pt]\beta _4 &=  \frac{C_1 P'}{\avrged{a^{-1}}^{2} \avrged{a}}, ~~~~~~ \beta _5 =\frac{2C_{10}+2C_{3}}{\avrged{a^{-1}} \avrged{a}}, ~~~~~~ \beta _6= \frac{ 2C_{1} C_{13}}{\avrged{a^{-1}}^{2} \avrged{a}^{2}}-\frac{C_{8}}{\avrged{a^{-1}}\avrged{a}^{2}}, \\[12pt] \beta _7 &=-\frac{2C_{1}C_{13}}{\avrged{a^{-1}}^{2} \avrged{a}^2}+\frac{4C_{7}}{\avrged{a^{-1}}^{2} \avrged{a} }-\frac{2C_{4}}{\avrged{a^{-1}}\avrged{a}} ,\\[12pt] \beta _{8} &=\frac{-3C_{5}+4C_{6}+C_{12}}{\avrged{a^{-1}}^{2}}+\frac{4C_{13}^2 + 4C_{13}C_{14}}{\avrged{a^{-1}}^{2} \avrged{a}}, \\[12pt] \beta _{9} &=\frac{2C_{7}P'}{ \avrged{a^{-1}}^{3}}+\frac{C_{1}P''}{\avrged{a^{-1}}^{2} \avrged{a}}-\frac{2 C_{1} C_{13} P'}{\avrged{a^{-1}}^{3} \avrged{a}},~~~~ \beta _{10}=\frac{C_{1}P''}{\avrged{a^{-1}}^{2} \avrged{a}}, ~~~~ \beta _{11}= \frac{C_{11}P'}{\avrged{a^{-1}}^{3} \avrged{a}}-\frac{C_{2} P'}{\avrged{a^{-1}}^{2} \avrged{a}}.
\end{align*}

\newpage
\bibliography{References}
\bibliographystyle{ieeetr}
\end{document}